\documentstyle[12pt,psfig,amssymb]{article}
\begin{document}
\title{Sylvester Waves in the Coxeter Groups}

\author{Leonid G. Fel${}^{\dag}$ and Boris Y. Rubinstein${}^{\ddag}$\\
\\${}^{\dag}$School of Physics and Astronomy, 
\\Raymond and Beverly Sackler Faculty
of Exact Sciences\\Tel Aviv University, Tel Aviv 69978, Israel\\
{\sl e-mail: lfel@ccsg.tau.ac.il}\\
and\\
${}^{\ddag}$Department of Engineering Sciences and Applied Mathematics\\            
       Northwestern University, 2145 Sheridan Road \\                               
       Evanston IL 60208-3125, U.S.A.}

\date{\today}

\maketitle

\def\be{\begin{equation}}
\def\ee{\end{equation}}
\def\p{\prime}

\begin{abstract}

A new recursive procedure of the calculation of partition numbers 
function $W(s,{\bf d}^m)$ 
is suggested. We find its zeroes and prove a lemma on the function 
parity properties. The explicit formulas of $W(s,{\bf d}^m)$ and 
their periods $\tau(G)$ for the irreducible Coxeter groups and a list for the first ten
symmetric group ${\cal S}_m$ are presented. 
A {\it least common multiple} ${\cal L}(m)$ of the series of the natural numbers 
1,2,..,$m$ plays a role of the period $\tau({\cal S}_m)$ of $W(s,{\bf d}^m)$
in ${\cal S}_m$. An asymptotic behaviour of 
${\cal L}(m)$ with $m \rightarrow \infty$ is found.

Pacs: Number theory, Invariant theory

\end{abstract}

\vskip 8 cm

\newpage
\noindent

\section{Introduction}

More than hundred years ago J.J.Sylvester stated \cite{sylv55,sylv59} and proved 
\cite{sylv82} a theorem about restricted partition number $W(s,{\bf d}^m)$ of 
positive integer $s$ with respect to the $m$-tuple of 
positive  integers ${\bf d}^m=\{d_1 , d_2, ... , d_m\}\;$:

{\bf Theorem.} {\it 
The number $W(s,{\bf d}^m)$ of ways in which $s$ can be composed
of (not necessarily distinct) $m$ integers $d_1,d_2,...,d_m$ is made up 
of a finite number of \underline{waves}
\begin{equation}
W(s,{\bf d}^m)=\sum_q^{max\;q} W_q(s,{\bf d}^m)\;\;,\;\;W_q(s,{\bf d}^m)=
\sum_k^{max\;k} W_{p_k|q}(s,{\bf d}^m)\;,
\label{syl1}
\end{equation}
where $q$ run over all distinct factors in $d_1,d_2,...,d_m$ and 
$W_{p_k|q}(s,{\bf d}^m)$ denotes the coefficient of $t^{-1}$ in the series 
expansion in ascending powers of $t$ of
\begin{equation}
F(s,{\bf d}^m,k;t)=\;e^{s w_k}\;\;\prod_{r=1}^m \frac{1}
{1-e^{d_r u_k}}\;\;,\;\;w_k=2\pi i\;
\frac{p_k}{q}+t\;\;,\;\;u_k=2\pi i\; \frac{p_k}{q}-t\;,
\label{syl2}
\end{equation}
and $p_1,p_2,...,p_{max\;k}$ are all numbers (unity included) less than 
$q$ and prime to it.
}

$W(s,{\bf d}^m)$ is also a number of sets of positive integer solutions 
$(x_1,x_2,...,x_m)$ of equation $\sum_r^m d_r x_r =s$. It is known that $W(s,{\bf d}^m)$ 
is equal to the coefficient of $t^s$ in the expansion of generating function
\begin{equation}
M({\bf d}^m,t)=\prod_{r=1}^m\frac{1}{1-t^{d_{r}}}=\sum_{s=0}^{\infty}
W(s,{\bf d}^m)\;t^s \;.
\label{syl3}
\end{equation}
If the exponents $d_1,d_2,...,d_m$ become the series of integers $1,2,3,...,m$, the 
number of waves is $m$ and $W(s,{\bf d}^m)$ of $s$ is usually referred to as a restricted 
partition number ${\cal P}_m(s)$  of $s$ into parts none of which exceeds $m$.

Another definition of $W(s,{\bf d}^m)$ comes from the polynomial invariant
of finite reflection groups. Let $M({\bf d}^m,t)$ is a Molien function 
of such a group $G$, 
$d_{r}$ are the degrees of basic invariants, and $m$ is the number of basic invariants 
\cite{humph90}. Then $W(s,{\bf d}^m)$ gives a number of algebraic 
independent polynomial invariants of the $s$-degree for group $G$.

Throughout his papers J.J.Sylvester gave different names for $W(s,{\bf d}^m)$ : {\it quotity}, 
{\it denumerant}, {\it quot-undulant} and {\it quot-additant}. Sometime after he discarded 
some of them. Because of a wide usage of $W(s,{\bf d}^m)$ not only as a partition number we 
shall call $W(s,{\bf d}^m)$ a {\it Sylvester wave}.

The Sylvester theorem is a very powerful tool not only in the trivial situation 
when $m$ is finite but also it was used for the purposes of asymptotic evaluations 
${\cal P}_m(s)$, as well as for the main term of the Hardy-Ramanujan formulas for
unrestricted partition number ${\cal P}(s)$ \cite{szeker51}.

Recent progress in the self-dual problem of effective isotropic conductivity in
two-dimensional three-component regular checkerboards \cite{fel20} and its further
extension on the $m$-component anisotropic cases \cite{fel21} have shown an existence
of algebraic equations with permutation invariance with respect to the action of
the finite group $G$ permuting $m$ components. $G$ is a subgroup of symmetric
group ${\cal S}_m$ and the coefficients in the equations are build out of algebraic
independent polynomial invariants for group $G$. Here $W(s,{\bf d}^m)$ measures 
a degree of non-universality of the algebraic solution with respect to the
different kinds of $m$-color plane groups.

Several proofs of Sylvester theorem are known \cite{sylv82},\cite{dicks50}. All of them make
use of the 
Cauchy$^,$s theory of residues. The recursion relations imposed on $W(s,{\bf d}^m)$ provide 
a combinatorial version of Sylvester formula. The classical example for the elementary 
(complex-variable-free) derivation was shown by Erd$\ddot{o}$s \cite{erd42} for the main 
term of the Hardy-Ramanujan formula. Recently an elementary derivation of Szekeres$^,$ 
formula for $W(s,{\bf d}^m)$ based on the recursion satisfied by $W(s,{\bf d}^m)$ was 
elaborated in \cite{canf96}. In this paper we give a new derivation of the Sylvester waves 
based on the recursion relation for $W(s,{\bf d}^m)$. We find also its {\it zeroes} 
and prove a lemma on parity properties of the Sylvester waves. Finally we present 
a list of the first ten Sylvester waves $W(s,{\cal S}_m)\;,m=1,...,10$ for symmetric groups
${\cal S}_m$ and for all Coxeter groups. In the Appendix we prove a
conjecture on asymptotic behaviour
of the {\it least common multiple} {\sf lcm}($1,2,...,N$) of the series of 
natural numbers.

\section{Recursion relation for $W(s,{\bf d}^m)$.}
We start with a recursion that follows from (\ref{syl3})
\begin{equation}
M({\bf d}^m,t)-M({\bf d}^{m-1},t)=t^{d_m} M({\bf d}^m,t),
\label{rec1}
\end{equation}
and after inserting the series expansions into the last equation we arrive at
\begin{equation}
W(s,{\bf d}^m)=W(s,{\bf d}^{m-1})+W(s-d_m,{\bf d}^m)\;,\;\;d_m \leq s\;,
\label{rec2}
\end{equation}
where $s$ is assumed to be real. We apply now the recursive procedure (\ref{rec2})
several times
\begin{equation}
W(s,{\bf d}^m)=\sum_{p=0}^{r_m} W(s-p\cdot d_m,{\bf d}^{m-1})+W(s-(r_m+1)\cdot d_m,{\bf d}^{m})\;.
\label{rec3}
\end{equation}
Let us consider {\it the generic form} of $W(k\cdot\tau \{{\bf d}^m\}+s,{\bf d}^m),
\;s < \tau \{{\bf d}^m\}$ where $k,\;s$ and $\tau \{{\bf d}^m\}$ are the independent  positive integers.
We will choose them in such a way that
\begin{equation}
k\cdot\tau \{{\bf d}^m\}+s-(r_m+1)\cdot d_m=(k-1)\cdot\tau \{{\bf d}^m\}+s\;,\;\;\;\; 
\Rightarrow \;\;\;\; \tau \{{\bf d}^m\}=(r_m+1)\cdot d_m\;.
\label{rec4}
\end{equation}
Thus the relation (\ref{rec3}) reads
\begin{eqnarray}
W(k\cdot\tau \{{\bf d}^m\}+s,{\bf d}^m)&=&W(\;(k-1)\cdot\tau \{{\bf d}^m\}+s,{\bf d}^m)+\nonumber \\
&&\sum_{p=0}^{\delta_m-1} W(k\cdot\tau \{{\bf d}^m\}-p\cdot d_m+s,{\bf d}^{m-1})\;\;,\;\;\;
\delta_m=\frac{\tau \{{\bf d}^m\}}{d_m}\;.
\label{rec5}
\end{eqnarray}
As follows from (\ref{rec4}), in order to return via the recursive procedure from
$W(k\cdot\tau \{{\bf d}^m\}+s,{\bf d}^m)$ to $W(\;(k-1)\cdot\tau \{{\bf d}^m\}+s,{\bf d}^m)$ we must use
$\tau \{{\bf d}^m\}$ which have $d_m$ as a divisor. Due to the arbitrariness of $d_m$ it is easy to
conclude that all exponents $d_1,d_2,...,d_m$ serve as the divisors of $\tau \{{\bf d}^m\}$. In other
words $\tau \{{\bf d}^m\}$ is the {\it least common multiple} {\sf lcm} of the exponents
$d_1,d_2,...,d_m$
\begin{equation}
\tau \{{\bf d}^m\}={\sf lcm}(d_1,d_2,...,d_m)\;.
\label{rec6}
\end{equation}
Actually $\tau \{{\bf d}^m\}$ does play a role of the "{\it period }" of $W(s,{\bf d}^m)$. 
But strictly speaking it is not a periodic function with respect to the integer variable $s$ as could 
be seen from (\ref{rec5}). The rest of the paper clarifies this hidden periodicity.

As we have mentioned above, $W(s,{\bf d}^m)$ gives a number of algebraic independent polynomial 
invariants of the $s$-degree for the group $G$. The situation becomes more transparent if we deal with
the irreducible Coxeter group where the degrees $d_{r}$ and the number of basic invariants $m$ are well
known.
\begin{center}
{\bf Table$\;$1.} The "{\it periods }" $\tau(G)$ of $W(s,{\bf d}^m)$  for
the irreducible Coxeter groups.
\begin{tabular}{|c||c|c|c|c|c|c|c|c|c|c|c|c|} \hline\hline
$G\;$ & $A_m$ & $B_m$ & $D_m$ & $G_2$ & $F_4$ & $E_6$ \\ \hline
$\tau(G)$ &  ${\cal L}(m+1)$ & 2${\cal L}(m)$ & 2${\cal L}(m)$ & $\;\;\;\;6\;\;\;\;$ & 24 & 360 \\
\hline\hline 
$G\;$ & $E_7$ & $E_8$ & $H_3$ & $H_4$ & $I_2$(2m) & $I_2$(2m+1)\\ \hline
$\tau(G)$ & 2520 & 2520 & 30 & 60 & 2$\;$m & 2$\;$(2m+1) \\ \hline\hline
\end{tabular}
\end{center}
where ${\cal L}(m)$={\sf lcm}(1,2,3,...,$m$) is the {\it least common multiple} of the series of
the natural numbers. 

${\cal L}(m)$ can be viewed as $\tau({\cal S}_m)$ for symmetric group ${\cal S}_m$ or, 
in other words, as a "{\it period }" of the restricted partition number ${\cal P}_m(s)$.
This makes it possible to pose a question about asymptotic behaviour of 
$\tau({\cal S}_m)$ with $m \rightarrow \infty$.
${\cal L}(m)$ is a very fast growing function: ${\cal L}(1)$=1, ${\cal L}(10)$=2520, 
${\cal L}(20)$=232792560, ${\cal L}(30)$=2329089562800 etc. Actually $\frac{\ln {\cal L}(m)}{m}$ 
oscillates infinitely many times around 1 and the function ${\cal L}(m)$ has an exponential increase 
with the asymptotic law 
\footnote{It seems to be strange but we have not found throughout the textbooks on number
theory any discussion about the asymptotics of {\sf lcm}(1,2,3,...,$m$). The formula
(\ref{rudnik1}) was conjectured by one of the authors (LGF) based on the numerical 
calculations and proved by Z.Rudnick (Tel-Aviv Univ., Israel) which had communicated this proof
to us. The proof is given in Appendix A.}
\begin{equation}
\lim_{m\to\infty} \frac{\ln {\cal L}(m)}{m}= 1.
\label{rudnik1}
\end{equation}

\section{Polynomial representation for $W(s,{\bf d}^m)$.}

Making use of the relations (\ref{rec5},\ref{rec6}) we obtain the exact formula for
$W(k\cdot\tau \{{\bf d}^m\}+s,{\bf d}^m)$ for different ${\bf d}^m$. We will treat it in an
ascending order in the number $m$ of exponents. The first steps are simple and they yield
$\\ \underline{{\bf d}^1=(d_1)}$ $\;,\;\;\tau\{{\bf d}^1\} > s\geq 0$
\begin{equation}
W(k\cdot d_1+s,{\bf d}^1)=W(s,{\bf d}^1)=\Psi_{d_1}(s)=
  \left\{ \begin{array}{ll}
         1\;\;,   & \mbox{$s=0 \pmod{d_1}$}  \\
         0\;\;,   & \mbox{$s\neq 0 \pmod{d_1}$}
                           \end{array}\right.
\label{d1}
\end{equation}
$\Psi_{d_1}(s)$ may be represented as a sum of prime roots of unit of degree $d_1$:
$$
\Psi_{d_1}(s) =
\frac{1}{d_1}\sum_{k=0}^{d_1-1} \exp(\frac{2 \pi i k s}{d_1}) =
\frac{1}{d_1}
\left\{
\begin{array}{ll}
1+\cos \pi s+2\sum_{k=1}^{d_1/2-1} \cos \frac{2 \pi k s}{d_1}, &
\mbox{even } d_1 \\
1+2\sum_{k=1}^{(d_1-1)/2} \cos \frac{2 \pi k s}{d_1}, & \mbox{odd } d_1
\end{array}
\right.
.
$$
$\\ \underline{{\bf d}^2=(d_1,d_2)}$ $\;,\;\;\tau\{{\bf d}^2\} > s\geq 0$
\begin{equation}
W(k\cdot \tau\{{\bf d}^2\}+s,{\bf d}^2)=W(s,{\bf d}^2)+k\cdot
\sum_{p=0}^{\delta_2-1} W(|s-p\;d_2|,{\bf d}^1)\;\;.
\label{d2}
\end{equation}
$\underline{{\bf d}^3=(d_1,d_2,d_3)}$ $\;,\;\;\tau\{{\bf d}^3\} > s\geq 0$ 
\begin{eqnarray}
W(k\cdot \tau\{{\bf d}^3\}+s,{\bf d}^3)&=&W(s,{\bf d}^3)+k\cdot \sum_{p=0}^{\delta_3-1}
W(|s-p\;d_3|,{\bf d}^2)+\nonumber \\ 
&&\frac{k(k+1)}{2}\;\frac{\tau\{{\bf d}^3\}}{\tau\{{\bf
d}^2\}}\sum_{p=0}^{\delta_3-1} \sum_{q=0}^{\delta_2-1} W(|s-p\;d_3-q\;d_2|,{\bf d}^1)\;. 
\label{d3}
\end{eqnarray} 
Now it is simple to deduce by induction that in the general case $W(k\cdot
\tau\{{\bf d}^{m}\}+s,{\bf d}^{m})$ has a polynomial representation with respect to $k$
\begin{equation} 
W(k\cdot \tau\{{\bf d}^{m}\}+s,{\bf d}^{m})=A_{m-1}^{m}(s)\;k^{m-1}\;+\;A_{m-2}^{m}(s)\;
k^{m-2}\;+ ...+A_{1}^{m}(s)\;k\;+\;A_{0}^{m}(s,{\bf d}^{m})\;, 
\label{dgen1} 
\end{equation} 
where $A_{m-r}^{m}(s)$ is based on the ${\tau\{{\bf d}^r\}}$-periodic functions as well as the 
entire $W(s,{\bf d}^{m})$ is based on the ${\tau\{{\bf d}^m\}}$-periodic functions. 
The coefficient in the leading term can be written in a closed form
\begin{eqnarray}
A_{m-1}^m (s)=\frac{1}{(m-1)!}\cdot \frac{\tau^{m-2}\{{\bf d}^m\}}{\tau\{{\bf d}^2\}\cdot
\tau\{{\bf d}^3\}\cdot ... \cdot \tau\{{\bf d}^{m-1}\}}\times \nonumber \\
\;\;\;\;\;\;\;\;\;\times\sum_{p=0}^{\delta_m-1}\sum_{q=0}^{\delta_{m-1}-1}\;...\;
\sum_{v=0}^{\delta_2-1}W(|s-p\;d_m-q\;d_{m-1}-...-v\;d_2|,{\bf d}^1)\;.
\label{dgen2}
\end{eqnarray}
With $d_1=1$ we have $W(|s-p\;d_m-q\;d_{m-1}-...-v\;d_2|,1)=1\;$,
which makes $A_{m-1}^m (s)$ independent of $s$ and gives an asymptotics 
of $W(s,{\bf d}^{m})$ for $s \gg m$
\begin{equation}
A_{m-1}^m (s)=\frac{\tau^{m-1}\{{\bf d}^m\}}{(m-1)!\;m!}\;,\;\;\;\;
W(s,{\bf d}^{m}) \stackrel{s \rightarrow \infty} \simeq \frac{s^{m-1}}{(m-1)!\;m!}\;.
\label{dgen3}
\end{equation}
Now we are ready to prove the statement about splitting of $W(s,{\bf d}^{m})$
into periodic and non-periodic parts. 

{\bf Lemma 3.1.}
{\it
The Sylvester wave $W(s,{\bf d}^{m})$ can be represented in the following way
\begin{equation}
W(s,{\bf d}^{m})=Q_m^m(s)+\sum_{j=1}^{m-1} Q_{j}^m(s)\cdot s^{m-j}\;,
\label{syl7}
\end{equation}
where $Q_j^m(s)$ is a periodic function with the period $\tau\{{\bf d}^{j}\}=
{\sf lcm}(d_1,d_2,...,d_j)$.
}
\vskip .5 cm
\underline{{\it Proof.}}$\;\;\;\;$
We start with the identity for the polynomial
representation for $W(k\cdot \tau\{{\bf d}^{m}\}+s,{\bf d}^{m})$
$$
W((k+1)\cdot \tau\{{\bf d}^{m}\}+s,{\bf d}^{m})=W(k\cdot \tau\{{\bf d}^{m}\}+s+
\tau\{{\bf d}^{m}\},{\bf d}^{m})\;,
$$
that can be transformed, using (\ref{dgen1}), into
\begin{eqnarray}
A_{m-1}^{m}(s)\;(k+1)^{m-1}\;+\;A_{m-2}^{m}(s)\;(k+1)^{m-2}\;+
...+A_{1}^{m}(s)\;(k+1)\;+\;W(s,{\bf d}^{m})=\nonumber \\
A_{m-1}^{m}(s+\tau\{{\bf d}^{m}\})\;k^{m-1}\;+\;
A_{m-2}^{m}(s+\tau\{{\bf d}^{m}\})\;k^{m-2}\;+...+
A_{1}^{m}(s+\tau\{{\bf d}^{m}\})\;k\;+\nonumber \\
W(s+\tau\{{\bf d}^{m}\},{\bf d}^{m})\;.
\label{syl2a}
\end{eqnarray}
The last identity generates a finite number of coupled difference equations for 
the coefficients $A_{r}^{m}(s)$
\begin{equation}
A_{m-r}^m(s+\tau\{{\bf d}^{m}\})=\sum_{j=1}^{r}C_{m-j}^{m-r}\cdot 
A_{m-j}^m(s)\;,\;\;\;1 \leq r \leq m\;,
\label{syl30}
\end{equation}
where $C^{k}_{n}$  denotes a binomial coefficient. The first equation $(r=1)$
\[ A_{m-1}^{m}(s+\tau\{{\bf d}^{m}\})=A_{m-1}^{m}(s)\;\]
declares that $A_{m-1}^{m}(s)$ is an arbitrary $\tau\{{\bf d}^{m}\}$-periodic function.
We can specify the last statement taking into account (\ref{dgen1}) that actually 
$A_{m-1}^{m}(s)$ is $\tau\{{\bf d}^1\}$-periodic function which will be denoted as 
$Q_1^m(s)$. The second equation $(r=2)$ 
\[A_{m-2}^{m}(s+\tau\{{\bf d}^{m}\})=A_{m-2}^{m}(s)+(m-1)\cdot A_{m-1}^{m}(s)\;\]
can be solved completely
\begin{equation}
A_{m-2}^{m}(s)=Q_2^m(s)+(m-1)\cdot s \cdot Q_1^m(s)\;,
\label{syl4}
\end{equation}
where $Q_2^m(s+\tau\{{\bf d}^{2}\})=Q_2^m(s)$. Continuing this procedure, it is not
difficult to prove by induction that for any $r$ we have
\begin{equation}
A_{m-r}^{m}(s)=\sum_{j=1}^r C_{m-j}^{m-r}\cdot 
Q_j^m(s)\cdot s^{r-j}\;,
\label{syl6}
\end{equation}
where $Q_j^m(s+\tau\{{\bf d}^{j}\})=Q_j^m(s)$. 
Since $W(s,{\bf d}^{m})=A_0^{m}(s)$ we arrive finally at (\ref{syl7}) 
by inserting $r=m$ into equation (\ref{syl6}), that splits 
$W(s,{\bf d}^{m})$, in accordance with the Sylvester theorem, into periodic and
non-periodic parts.$\;\;\;\;\;\;\blacksquare$

\section{Partition identities and zeroes of $W(s,{\bf d}^{m})$.}

In this section we assume that the variable $s$ has only integer values.

Consider a new quantity 
\begin{equation}
V(s,{\bf d}^m) = W(s-\xi \{{\bf d}^m\},{\bf d}^m)\;,
\;\;\xi \{{\bf d}^m\}=\frac{1}{2} \sum_{i=1}^m d_i\;.
\label{transform1}
\end{equation}

{\bf Lemma 4.1.}
{\it 
$V(s,{\bf d}^m)$ has the following parity properties:
\begin{equation}
V(s,{\bf d}^{2m})=-V(-s,{\bf d}^{2m}), \ \ \
V(s,{\bf d}^{2m+1})=V(-s,{\bf d}^{2m+1}).
\label{lem11}
\end{equation}}

\underline{{\it Proof.}}$\;\;\;\;$
A basic recursion relation (\ref{rec2}) can be rewritten for $V(s,{\bf d}^m)$
\begin{equation}
V(s,{\bf d}^m) - V(s-d_m,{\bf d}^m) = V(s-\frac{d_m}{2},{\bf d}^{m-1})\;.
\label{1}
\end{equation}
The last relation produces two equations in a new variable $q=s-\frac{d_m}{2}$
\begin{eqnarray}
V(q,{\bf d}^{m-1}) &=& V(q+\frac{d_m}{2},{\bf d}^m) - V(q-\frac{d_m}{2},{\bf d}^{m}), \cr
V(-q,{\bf d}^{m-1}) &=& V(-q+\frac{d_m}{2},{\bf d}^m) - V(-q-\frac{d_m}{2},{\bf d}^{m})\;.
\label{2}
\end{eqnarray}
Hence if $V(q,{\bf d}^{m})$ is an even function of $q$, then $V(q,{\bf d}^{m-1})$ is an odd one, 
and vice versa. Because $V(q,{\bf d}^{1})$ is an even function, 
we arrive at (\ref{lem11}). $\;\;\;\;\;\;\blacksquare$

\vspace{.3 cm}
{\bf Corollary.} {\it
If $\;\;\;s_1+s_2+2 \xi \{{\bf d}^m\}=0\;\;\;$, then
$$
W(s_1,{\bf d}^m)=(-1)^{m+1} W(s_2,{\bf d}^m)
$$
}

\underline{{\it Proof.}}$\;\;\;\;$
This follows from the parity properties and after substitution two new variables
$s_1=s-\xi \{{\bf d}^m\}\;,\;\;s_2=-s-\xi \{{\bf d}^m\}\;$ into
(\ref{lem11}). $\;\;\;\;\;\;\blacksquare$

\vspace{.3 cm}
{\bf Lemma 4.2.} {\it 
Let $m$-tuple $\{{\bf d}^m\}$ generates the Sylvester wave $W(s,{\bf d}^{m})$ . 
Then for every integer $p$ a $m$-tuple $\{p\cdot {\bf d}^m\}=\{p d_1,p d_2,...,p d_m\}$
generates the following Sylvester wave}
\begin{equation}
W(s,p\cdot {\bf d}^m)=\Psi_p(s)\cdot W(\frac{s}{p},{\bf d}^{m})\;, \;\;\mbox{or}\;\;\;\;
V(s,p\cdot {\bf d}^m)=\Psi_p(s-p \xi \{{\bf d}^m\})\cdot V(\frac{s}{p},{\bf d}^{m})\;,
\label{lem2}
\end{equation}
{\it where the periodic function $\Psi_p(s)=\Psi_p(s+p)$ is defined in (\ref{d1}).
}

\vskip .6 cm
\underline{{\it Proof.}}$\;\;\;\;$
According to the definition (\ref{syl3})
\[ \sum_s W(s,p\cdot {\bf d}^m)\cdot t^s=\sum_s W(s,{\bf d}^m)\cdot t^{ps}=
\sum_{s^{\prime}} W(\frac{s^{\prime}}{p},{\bf d}^m)\cdot t^{s^{\prime}} \]
Equating powers of $t$ in the latter equation and taking into account that $s^{\prime}/p$
must be integer we obtain (\ref{lem2}). $\;\;\;\;\;\;\blacksquare$

\vspace{.3 cm}
{\bf Lemma 4.3.} {\it
Let $m$-tuple $\{{\bf d}^m\}$ generates the Sylvester wave $W(s,{\bf d}^{m})$ .
Then $W(s,{\bf d}^{m})$ has the following zeroes:}
\begin{itemize}
\item {\it If all exponents $d_r$ are mutually prime numbers, then the zeroes 
${\mathfrak{s}}_0({\bf d}^{m})$ read}
\begin{eqnarray}
{\mathfrak{s}}_0({\bf d}^{m})&=&-1,-2,...,-\sum_{r=1}^m d_r+1\;,\;\;\;\;\mbox{if}\;\;\; m
=2k+1\;,\nonumber \\
{\mathfrak{s}}_0({\bf d}^{m})&=&-1,-2,...,-\sum_{r=1}^m d_r+1\;,\; -\xi \{{\bf d}^m\}\;,\;\;\;\;
\mbox{if}\;\;\; m =2k\;;
\label{lem31}
\end{eqnarray}
\item {\it If all exponents $d_r$ have a maximal common factor $p$, then $W(s,{\bf d}^{m})$ 
has infinite number of zeroes ${\mathfrak{S}}_1({\bf d}^{m})$ which are distributed in the 
following way} 
\begin{equation}
{\mathfrak{S}}_1({\bf d}^{m})={\mathfrak{s}}_1({\bf d}^{m}) 
\cup \{ {\bf \Bbb{Z}} / p {\bf \Bbb{ Z}} \} \;\;\;,\nonumber \\
\label{lem32}
\end{equation}
{\it where $\{ {\bf \Bbb{Z}} /p {\bf \Bbb{ Z}} \}$ denotes a set of integers ${\bf \Bbb{Z}}$ with 
deleted integers of modulo} $p$ 
\begin{equation}
\{ {\bf \Bbb{Z}} / p {\bf \Bbb{ Z}} \}=\{...,-p-1,-p+1,...,-1,1,...,p-1,p+1,...\} \nonumber
\end{equation}
{\it and} 
\begin{eqnarray}
{\mathfrak{s}}_1({\bf d}^{m})&=& -p,-2p,...,-\sum_{r=1}^m d_r+p\;,\;\;\;\;\mbox{if}\;\;\; 
 m =2k+1\;,\nonumber \\
{\mathfrak{s}}_1({\bf d}^{m})&=& -p,-2p,...,-\sum_{r=1}^m d_r+p,\;- \xi \{{\bf d}^m\}\;,\;\;
\mbox{if}\;\;\;\; m = 2k\;.
\label{lem33}
\end{eqnarray}
\end{itemize}

\underline{{\it Proof.}}$\;\;\;\;$
Consider again the relation (\ref{rec3}) which we rewrite as follows
\begin{equation}
\sum_{s=0}^{\infty}W(s,{\bf d}^m)\cdot t^s=
\frac{1}{1-t^{d_m}}\cdot \sum_{s'=0}^{\infty}W(s',{\bf d}^{m-1})\cdot t^{s'}
\label{rec3zer}
\end{equation}
assuming that the exponents in ${\bf d}^m$ 
are sorted in the ascending order. Note that the influence of the new $d_m$ exponent 
appears only in terms $t^s$ with $s \ge d_m$. This enables us to deduce that the values of 
$W(s,{\bf d}^{m-1})$ and $W(s,{\bf d}^{m})$ coincide at 
integer positive values $s=0,1,\ldots,d_m-1$. This means that for $0 \le s \le d_m-1$ we have
$W(s,{\bf d}^{m})=W(s,{\bf d}^{m-1})$. Recalling the main recursion relation (\ref{rec2}) 
we conclude that $$
W(s,{\bf d}^{m}) = 0 \ \ (-d_m \le s \le -1).
$$
Using the last relation for $m$ and $m-1$ in (\ref{rec2}) we can find also
$$
W(s-d_m,{\bf d}^{m}) = 0 \ \ (-d_{m-1} \le s \le -1) \ \  \Rightarrow \ \
W(s,{\bf d}^{m}) = 0 \ \ (-d_{m-1}-d_m \le s \le -1).
$$
Repeating this procedure and taking into account that at the last step it 
leads to the zeroes of $\Psi_{d_{1}}$ which are located at $(1-d_{1} \le s \le -1)$, 
we get the set of the zeroes for $W(s,{\bf d}^{m})$ with odd number of exponents $m = 2k+1$
\begin{equation}
W(s,{\bf d}^{m}) = 0 \ \ (1-\sum_{i=1}^m d_i \le s \le -1).
\label{zerfin}
\end{equation}
The eveness of $m$ gives one more zero of $W(s,{\bf d}^{m})$ which arises from the parity properties
of $V(s,{\bf d}^{m})$, namely, $V(0,{\bf d}^{2k})=0$. The last equality immediately generates a zero 
$-\xi \{{\bf d}^{2k}\}$ of $W(s,{\bf d}^{2k})$ that together with (\ref{zerfin}) proves the first
part (\ref{lem31}) of Lemma 3.

The second part of Lemma 3 follows from (\ref{lem2}) and from the first part of (\ref{lem31})
 because a set of integers $\{ {\bf \Bbb{Z}} /p {\bf \Bbb{ Z}} \}$ represents the zeroes
of the periodic function $\Psi_p(s)$. $\;\;\;\;\;\;\;\;\;\;\blacksquare$

The complexity of the exponents sequence $\{{\bf d}^m\}$ and its large length make the calculative
procedure of restoration of $Q_{j}^m(s)$ very cumbersome. Therefore it is important to find the inner
properties of $\{{\bf d}^m\}$ when this procedure could be essentially reduced.

\vskip .3 cm
{\bf Lemma 4.4.} {\it Let $m$-tuple $\{{\bf d}^m\}=\{d_1,d_2,...,d_r,d_r,...,d_m\}$
contains an exponent $d_r$ twice. Then the Sylvester wave $V(s,{\bf d}^{m})$ is related
to the Sylvester wave $V(s,{\bf d}^{m_1})$ produced by the the non-degenerated 
tuple $\{{\bf d}^{m_1}\}=\{d_1,d_2,...,d_r,..,.d_m,2d_r\}$ as follows}
\begin{equation}
V(s,{\bf d}^m)=V(s-\frac{d_r}{2},{\bf d}^{m_1})+
V(s+\frac{d_r}{2},{\bf d}^{m_1})\;.
\label{lem4}
\end{equation}

\vskip .5 cm
\underline{{\it Proof.}}$\;\;\;\;$
According to the definition (\ref{syl3})
$$
(1+t^{d_r})\cdot \sum_s W(s,{\bf d}^{m_1})\cdot t^s=
\sum_s W(s,{\bf d}^m)\cdot t^s\;.
$$
Taking into account that $\xi \{{\bf d}^{m_1}\}-\xi \{{\bf d}^m\}=d_r/2$ and equating powers 
of $t$ in the latter equation we obtain the stated relation (\ref{lem4}) according to the definition
(\ref{transform1}). $\;\;\;\;\;\;\blacksquare$

We will make worth of relation (\ref{lem4}) during the evaluation of the expression $V(s,{\bf d}^{m})$
for the  Coxeter group {\sl D}$_m$.

\section{Recursion formulas for $V(s,{\bf d}^{m})$.}
The shift (\ref{transform1}) transforms the relation (\ref{rec5}) into
\begin{equation}
V(s+\tau \{{\bf d}^m\},{\bf d}^{m})=V(s,{\bf d}^{m}) + 
\sum_{p=0}^{\delta_m-1} V(s+\tau \{{\bf d}^m\}-\lambda_p\cdot d_m,{\bf d}^{m-1})\;,\;\;
\lambda_p=p+\frac{1}{2}
\label{3}
\end{equation}
and the relation (\ref{syl7}) into
\begin{equation}
V(s,{\bf d}^{m})=R^{m}_m(s)+\sum_{j=1}^{m-1}R^{m}_j(s)\cdot s^{m-j}\;,
\label{31}
\end{equation}
where
\[ R^{m}_j(s)=\sum_{i=1}^j C^{j-i}_{m-i}\cdot(-\xi \{{\bf d}^m\})^{j-i}\cdot
Q_i^m(s-\xi \{{\bf d}^m\})\;,\]
i.e., $R^{m}_1(s)=Q^{m}_1(s-\xi \{{\bf d}^m\})\;;\;R^{m}_2(s)=Q^{m}_2(s-\xi \{{\bf d}^m\})-
(m-1)\cdot \xi \{{\bf d}^m\}\cdot  Q^{m}_1(s-\xi \{{\bf d}^m\})$ etc. This means that the functions 
$R^{m}_j(s)$ and $Q^{m}_j(s)$ have the same period $\tau \{{\bf d}^j\}$.

Inserting the expansion (\ref{31}) into the relation (\ref{3}) and equating powers of $s$ we can
obtain for $k=1,2,\ldots,m-1$
\begin{equation}
\sum_{j=1}^{k} C^{m-1-k}_{m-j}\cdot R^{m}_j(s)\cdot \tau \{{\bf d}^m\}^{k+1-j} =
\sum_{p=0}^{\delta_m-1} \sum_{j=1}^{k} R^{m-1}_j(s-\lambda_p\cdot d_m)\cdot 
C^{m-1-k}_{m-1-j}\cdot (\tau \{{\bf d}^m\}-\lambda_p\cdot d_m)^{k-j}.
\label{master1}
\end{equation}
For the first successive values of $k$ the latter equation (\ref{master1}) gives
\begin{eqnarray}
R^{m}_1(s) &=& \frac{1}{(m-1)\cdot \tau \{{\bf d}^m\}}\sum_{p=0}^{\delta_m-1} 
R^{m-1}_1(s-\lambda_p\cdot d_m)\;,\nonumber \\
R^{m}_2(s) &=& \frac{1}{(m-2)\cdot \tau \{{\bf d}^m\}}\sum_{p=0}^{\delta_m-1}
R^{m-1}_2(s-\lambda_p\cdot d_m)+
\sum_{p=0}^{\delta_m-1}(\frac{1}{2}-\frac{\lambda_p}{\delta_m})\cdot
R^{m-1}_1(s-\lambda_p\cdot d_m)\;,\nonumber \\
R^{m}_3(s) &=& \frac{1}{(m-3)\cdot \tau \{{\bf d}^m\}}\sum_{p=0}^{\delta_m-1}
R^{m-1}_3(s-\lambda_p\cdot d_m)+
\sum_{p=0}^{\delta_m-1}(\frac{1}{2}-\frac{\lambda_p}{\delta_m})\cdot 
R^{m-1}_2(s-\lambda_p\cdot d_m)+\nonumber \\
&&\frac{m-2}{2}\cdot \tau \{{\bf d}^m\}\sum_{p=0}^{\delta_m-1}
(\frac{1}{6}-\frac{\lambda_p}{\delta_m}+\frac{\lambda^2_p}{\delta^2_m})\cdot 
R^{m-1}_1(s-\lambda_p\cdot d_m)\;.
\label{master2}   
\end{eqnarray}
It is easy to see that in the summands of the latter formulas (\ref{master2}) there appear 
the Bernoulli polynomials ${\cal B}_i(1-\frac{\lambda_p}{\delta_m})$ : ${\cal B}_0(x)=1,\;
{\cal B}_1(x)=x-1/2,\;{\cal B}_2(x)=x^2-x+1/6,\;
{\cal B}_3(x)=x^3-3/2\; x^2+1/2\; x\;$, etc \cite{bat53}.
Continuing the evaluation of the general expression for $R^{m}_{j}(s),\;1<j<m$, 
we arrive at

{\bf Lemma 5.1.}
{\it
$R^{m}_{j}(s)$ for $1 \leq j<m$ is given by the formula
\begin{equation}
R^{m}_j(s) =\frac{1}{m-j}\cdot 
\sum_{l=0}^{j-1} (\tau \{{\bf d}^m\})^{l-1}\cdot C^l_{m-1-j+l} \sum_{p=0}^{\delta_m-1} 
{\cal B}_l (1-\frac{\lambda_p}{\delta_m})\cdot R^{m-1}_{j-l}(s-\lambda_p\cdot d_m)\;.
\label{knf}
\end{equation}
}

\underline{{\it Proof.}}$\;\;\;\;$
Before going to the proof we recall two identities  for the Bernoulli polynomials
\cite{bat53},\cite{prud92}
\begin{equation}
{\cal B}_{l}(x+y) -{\cal B}_{l}(x)= \sum_{j=1}^l C_l^j\cdot y^j\cdot {\cal B}_{l-j}(x)\;,
\;\;\;{\cal B}_{l} (1+x) - {\cal B}_{l}(x) = l x^{l-1}\;.
\label{knf01}
\end{equation}
Using the definition (\ref{31}) we check that formula (\ref{knf}) satisfies (\ref{3}).
\begin{eqnarray}
&&V(s,{\bf d}^{m})  =  R^{m}_m(s)+
\sum_{j=1}^{m-1} s^{j} \sum_{l=j}^{m-1} C_l^j
\frac{(\tau \{{\bf d}^m\})^{l-j-1}}{l}
\sum_{p=0}^{\delta_m-1}
{\cal B}_{l-j}(1-\frac{\lambda_p}{\delta_m})R^{m-1}_{m-l}(s-\lambda_p d_m) =
\nonumber \\
&&
R^{m}_m(s)+ \sum_{l=1}^{m-1} \frac{(\tau \{{\bf d}^m\})^{l-1}}{l}
\sum_{p=0}^{\delta_m-1} R^{m-1}_{m-l}(s-\lambda_p d_m)
\sum_{j=1}^{l} C_l^j
\left(
\frac{s}{\tau \{{\bf d}^m\}}
\right)^j
{\cal B}_{l-j}(1-\frac{\lambda_p}{\delta_m}) = \nonumber \\
&&
R^{m}_m(s)+ \sum_{l=1}^{m-1} \frac{(\tau \{{\bf d}^m\})^{l-1}}{l}
\sum_{p=0}^{\delta_m-1} R^{m-1}_{m-l}(s-\lambda_p d_m)
\left [
{\cal B}_{l}(1+\frac{s-\lambda_p d_m}{\tau \{{\bf d}^m\}}) -
{\cal B}_{l}(1-\frac{\lambda_p}{\delta_m})
\right],
\label{knf2}
\end{eqnarray}
where we use the first of the identities (\ref{knf01}). Having in mind the $\tau
\{{\bf d}^m\}$-periodicity of functions $R^{m}_{j}(s)$ and $R^{m-1}_{j}(s)$ and 
the second identity (\ref{knf01})
we may rewrite the difference in the l.h.s of relation (\ref{3}) in the following
form:
\begin{eqnarray}
&&V(s,{\bf d}^{m})-V(s-\tau \{{\bf d}^m\},{\bf d}^{m}) = \\
&&
\sum_{l=1}^{m-1} \frac{(\tau \{{\bf d}^m\})^{l-1}}{l}
\sum_{p=0}^{\delta_m-1} R^{m-1}_{m-l}(s-\lambda_p d_m)
\left [
{\cal B}_{l}(1-\frac{\lambda_p}{\delta_m}+\frac{s}{\tau \{{\bf d}^m\}}) -
{\cal B}_{l}(-\frac{\lambda_p}{\delta_m}+\frac{s}{\tau \{{\bf d}^m\}})
\right] \nonumber \\
&&
\sum_{l=1}^{m-1} \frac{(\tau \{{\bf d}^m\})^{l-1}}{l}
\sum_{p=0}^{\delta_m-1} R^{m-1}_{m-l}(s-\lambda_p d_m) l
\left(
\frac{s-\lambda_p d_m}{\tau \{{\bf d}^m\}}
\right)^{l-1} = \nonumber \\
&&
\sum_{p=0}^{\delta_m-1}
\sum_{l=0}^{m-2} (s-\lambda_p d_m)^{l} R^{m-1}_{m-1-l}(s-\lambda_p d_m) =
\sum_{p=0}^{\delta_m-1} V(s-\lambda_p d_m,{\bf d}^{m-1}).
\label{knf3} \nonumber\;\;\;\;\;\;\;\;\;\;\;\blacksquare
\end{eqnarray}

The formula (\ref{knf}) enables to restore all terms $R^m_k(s)$ except the last
$R^m_m(s)$. Actually we can learn about it from the following consideration.
Let us separate $R^m_{m-k}(s)$ in the following way
\begin{equation}
R^m_{m-k}(s)={\cal R}^m_{m-k}(s)+r^m_{m-k}(s)\;,\;\;0\leq k\leq m-1\;,
\label{knf6}
\end{equation}
where
\begin{eqnarray}
{\cal R}^m_{m-k}(s)&=&\sum_{l=1}^{m-k-1} \frac{(\tau \{{\bf d}^m\})^{l-1}}{l+k}
\cdot C_{l+k}^k \sum_{p=0}^{\delta_m-1}
{\cal B}_l (1-\frac{\lambda_p}{\delta_m})\cdot
R^{m-1}_{m-k-l}(s-\lambda_p\cdot d_m) \label{knf7}\\
r^m_{m-k}(s)&=&\frac{1}{k\cdot \tau \{{\bf d}^m\}} \sum_{p=0}^{\delta_m-1}
R^{m-1}_{m-k}(s-\lambda_p d_m)\;,\;r^m_{m-k}(s) = r^m_{m-k}(s - d_m)\;,\;
(k \neq 0) \label{knf7a}
\end{eqnarray}
The representation (\ref{knf6}) and $d_m$-periodicity of the function
$r^m_{m-k}(s)$ make possible to prove the following 

{\bf Lemma 5.2.}
{\it
$R^m_{m-k}(s)$ for $0\leq k\leq m-1$ and ${\cal R}^m_{m-k}(s)$ for $0 < k\leq m-1$
satisfy the recursion relation
\begin{eqnarray}
&&R^m_{m-k}(s)-R^m_{m-k}(s-d_m)={\cal R}^m_{m-k}(s)-{\cal R}^m_{m-k}(s-d_m)=\label{knf9}
\nonumber \\
&&\sum_{j=k+1}^{m-1} \left \{ (-d_m)^{j-k}\cdot C_j^k\cdot R^m_{m-j}(s-d_m)+
(-\frac{d_m}{2})^{j-1-k}\cdot C^k_{j-1}\cdot
R^{m-1}_{m-j}(s-\frac{d_m}{2})\right \}.
\end{eqnarray}}

\underline{{\it Proof.}}$\;\;\;\;$
Inserting (\ref{31}) into (\ref{1}), expanding the powers of binomials into sums and
equating the powers of $s$ in the latter equation we obtain the relation
(\ref{knf9}) for the function $R^m_{m-k}(s)$, $0\leq k\leq m-1$. Using the definition
(\ref{knf6}) we immediately arrive at the relation for the function ${\cal R}^m_{m-k}(s)$ ,
$0< k\leq m-1$.
$\;\;\;\;\;\;\blacksquare$

In the special case $k=0$ the general relation (\ref{knf9}) produces the recursion
for $R^m_m(s)$
\begin{equation}
R^{m}_m(s) - R^{m}_m(s-d_m)=
\sum_{j=1}^{m-1} \left \{(-d_m)^{j}\cdot R^{m}_{m-j}(s-d_m) +
(-\frac{d_m}{2})^{j-1}\cdot R^{m-1}_{m-j}(s-\frac{d_m}{2})\right \}.
\label{knf10}
\end{equation}
We can not use directly (\ref{knf7}) for $k=0$ since $r^m_m(s)$ can not be derived
from (\ref{knf7a}). But it is a good mathematical intuition to exploit the formula
(\ref{knf7}) for $k=0$ in order to prove

{\bf Lemma 5.3.}
{\it
${\cal R}^{m}_{m}(s)$ is given by the formula
\begin{equation}
{\cal R}^m_m(s)=
\sum_{l=1}^{m-1} \frac{(\tau \{{\bf d}^m\})^{l-1}}{l} \sum_{p=0}^{\delta_m-1}
{\cal B}_l (1-\frac{\lambda_p}{\delta_m})\cdot R^{m-1}_{m-l}(s-\lambda_p\cdot d_m)\;.
\label{knf11}
\end{equation}
}

\underline{{\it Proof.}}$\;\;\;\;$
In order to prove that ${\cal R}^m_m(s)$ given by (\ref{knf11})
satisfies the difference equation (\ref{knf10}) we
consider a difference ${\cal R}^m_m(s)-{\cal R}^m_m(s-d_m)=
\Delta_m(s) = \Delta_m^1(s) + \Delta_m^2(s)$:
$$
\Delta_m(s)=
\sum_{l=1}^{m-1} \frac{(\tau \{{\bf d}^m\})^{l-1}}{l}
\sum_{p=0}^{\delta_m-1}
{\cal B}_l (1-\frac{\lambda_p}{\delta_m})\cdot
\left [
R^{m-1}_{m-l}(s-\lambda_p d_m)-
R^{m-1}_{m-l}(s-\lambda_{p+1} d_m)
\right]
$$
with
$$
\Delta_m^1(s)=\sum_{l=1}^{m-1} \frac{(\tau \{{\bf d}^m\})^{l-1}}{l}
\left \{
{\cal B}_{l} (1-\frac{1}{2\delta_m})-
 {\cal B}_{l} (-\frac{1}{2\delta_m})\right\}\cdot
R^{m-1}_{m-l}(s-\frac{d_m}{2})\;,
$$
$$
\Delta_m^2(s)=\sum_{l=1}^{m-1} \frac{(\tau \{{\bf d}^m\})^{l-1}}{l}
\sum_{p=1}^{\delta_m} \left \{
{\cal B}_{l} (1-\frac{\lambda_p}{\delta_m})-
{\cal B}_{l} (1-\frac{\lambda_p}{\delta_m}+\frac{1}{\delta_m})
\right\}\cdot
R^{m-1}_{m-l}(s-\lambda_p d_m)\;.
$$
The first term $\Delta_m^1(s)$ is calculated with the help of one of the identities
(\ref{knf01}):
\begin{equation}
\Delta_m^1(s)=\sum_{l=1}^{m-1} (-\frac{d_m}{2})^{l-1}\cdot
R^{m-1}_{m-l}(s-\frac{d_m}{2})\;.
\label{knf11a}
\end{equation}
Using another identity from (\ref{knf01}) we may write for $\Delta_m^2(s)$:
$$
\Delta_m^2(s)=\sum_{l=1}^{m-1} \sum_{j=1}^l\frac{(\tau \{{\bf d}^m\})^{l-1}}{l}
\cdot C_l^j \cdot (-\frac{1}{\delta_m})^j \sum_{p=1}^{\delta_m}
{\cal B}_{l-j}(1-\frac{\lambda_{p-1}}{\delta_m})\cdot
R^{m-1}_{m-l}(s-\lambda_p d_m)\;.  
$$
Changing here summation order
$
\sum_{k=l+1}^{m-1} \sum_{j=l+1}^{k} = \sum_{j=l+1}^{m-1} \sum_{k=j}^{m-1}
$
and comparing the inner sum with (\ref{knf}) we arrive at
\begin{equation}
\Delta_m^2(s) =\sum_{j=1}^{m-1} (-d_m)^j\cdot R^{m}_{m-j}(s-d_m)
\label{knf15} 
\end{equation}
Then (\ref{knf11a}) and (\ref{knf15}) prove the Lemma.
$\;\;\;\;\;\;\blacksquare$

From this Lemma follows an existence of $d_m$-periodic function $r^m_m(s)=r^m_m(s-d_m)$ 
which could not be derived from (\ref{knf7a}). 
Unknown function $r^m_m(s)$ corresponds to vanishing  harmonics in the r.h.s. of equation 
(\ref{knf9}). We are free to choose any basic system of continuous 
$\tau \{{\bf d}^m\}$-periodic 
functions. This arbitrariness can affect behaviour of $W(s,{\bf d}^m)$ only for non-integer
$s$ that does not violate the recursion relation (\ref{rec2}). In the rest of the paper
we will choose a basic system of the simplest periodic functions {\sf sin} and {\sf cos}.

The function $r^m_m(s)$ corresponds to the harmonics of the type
$$
\left\{
\begin{array}{c}
\sin \\ \cos
\end{array}
\right\} \frac{2\pi n}{d_m}s
$$
Because the parity properties of $R^m_m(s)$ coincide with that of $V(s,{\bf d}^{m})$ itself
we can rewrite (\ref{31}) in the following form
\begin{eqnarray}
V(s,{\bf d}^{2m})&=&\sum_{j=1}^{2m-1}R^{2m}_j(s)\cdot s^{2m-j}+{\cal R}^{2m}_{2m}(s)+
\sum_{n} \rho^{2m}_n \cdot \sin \frac{2\pi n}{d_{2m}}s\;, \label{knf12} \\
V(s,{\bf d}^{2m+1})&=&\sum_{j=1}^{2m}R^{2m+1}_j(s)\cdot s^{2m+1-j}+{\cal R}^{2m+1}_{2m+1}(s)+
\sum_{n} \rho^{2m+1}_n \cdot \cos  \frac{2\pi n}{d_{2m+1}}s\;.\label{knf13}
\end{eqnarray}
In order to produce $r^m_m(s)$ we use some of zeroes $\mathfrak{s}$, described in the
preceding Section, constructing a system of linear equations 
for $[(m+1)/2]$ coefficients $\rho_n$;
$n$ runs from $1$ to $m/2$ in (\ref{knf12} and
from $0$ to $(m-1)/2$ in (\ref{knf13}). We use a trivial identity $V(\xi({\bf d}^{m}),{\bf
d}^{m})=1$, and choose the values of $s$ out of the set $\mathfrak{s}$, adding homogeneous
equations to arrive at a non-degenerate inhomogeneous system of linear equations. This
system is solved further to produce the final expression for corresponding
Sylvester wave. These explicit expressions are given in the next Section.
Appendix B presents two instructive examples of the above procedure.

\section{Sylvester waves $V(s,G)$.}

We start with the symmetric group ${\cal S}_m$ because of two reasons: first, 
of their relation with restricted partition numbers and , second, they arranged 
natural basis to utilize the Sylvester waves $V(s,G)$ in all Coxeter groups.

\subsection{Symmetric groups ${\cal S}_m$.}

Making use of the procedure developed in the previous section we present here
first ten Sylvester waves $V(s,{\cal S}_m)\;,m=1,...,10$.
\footnote{Having in mind the results of Sylvester \cite{sylv55},\cite{sylv59} and Glaisher 
\cite{glaish10} for restricted partition numbers for $m \le 9$  we repeat them
adding a formula for $m=10$.
The list of $V(s,{\cal S}_m)$ can be simply continued up to any finite $m$ with 
the help of the symbolic code written in {\it Mathematica} language
\cite{wolf96}.}

$\underline {G={\cal S}_m}$ $\;,\;\;\;d_r=1,2,3,..., m\;
,\;\;\;\xi({\cal S}_m)=\frac{m(m+1)}{4}\;,$
\begin{eqnarray}
V(s,{\cal S}_1) & = & 1, \nonumber \\
V(s,{\cal S}_2) & = & \frac{s}{2} - \frac{1}{4} \sin \pi s, \nonumber \\
V(s,{\cal S}_3) & = & \frac{s^2}{12} - \frac{7}{72} 
-\frac{1}{8} \cos \pi s +\frac{2}{9} \cos \frac{2 \pi s}{3}, \nonumber \\
V(s,{\cal S}_4) & = & \frac{s^3}{144} - \frac{s}{96}\cdot (5+3 \cos \pi s)
+\frac{1}{8} \sin \frac{\pi s}{2} 
-\frac{2}{9\sqrt{3}} \sin \frac{2 \pi s}{3}, \label{sym4} \\ 
V(s,{\cal S}_5) & = & \frac{s^4}{2880} - \frac{11\cdot s^2}{1152} 
- \frac{s}{64}\cdot  \sin \pi s+\frac{17083}{691200}  -
\frac{2}{27} \cos \frac{2 \pi s}{3} + \nonumber \\
&& \frac{1}{8\sqrt{2}} \cos \frac{\pi s}{2} +
\frac{2}{25} (-\cos \frac{2 \pi s}{5} + \cos \frac{4 \pi s}{5}), \nonumber \\
V(s,{\cal S}_6) & = & \frac{s^5}{86400} -
\frac{91\cdot s^3}{103680} + \frac{s^2}{768}\cdot \sin \pi s
+ \frac{s}{829440}\cdot (9191-10240 \cos \frac{2 \pi s}{3})- \nonumber \\
&&\frac{161}{9216} \sin \pi s -
\frac{1}{16\sqrt{2}} \sin \frac{\pi s}{2} 
-\frac{1}{81\sqrt{3}} \sin \frac{2 \pi s}{3} 
-\frac{1}{18} \sin \frac{\pi s}{3} - \nonumber \\
&&\frac{2}{25\sqrt{5}}(\sin \frac{\pi}{5} \sin \frac{4 \pi s}{5} +
\sin \frac{2 \pi}{5} \sin \frac{2 \pi s}{5}),\nonumber \\
V(s,{\cal S}_7) & = & \frac{s^6}{3628800} -
\frac{s^4}{20736} + \frac{s^2}{38400}\cdot (71+25\cos \pi s)
-\frac{s}{81\sqrt{3}}\cdot \sin\frac{2 \pi s}{3}- \nonumber \\
&&\frac{52705}{6096384}-\frac{77}{4608} \cos \pi s -
\frac{1}{32} \cos \frac{\pi s}{2} 
-\frac{5}{486} \cos \frac{2 \pi s}{3} 
-\frac{1}{18} \cos \frac{\pi s}{3} + \nonumber \\
&&\frac{2}{25\sqrt{5}}
(\cos \frac{2 \pi s}{5}-\cos \frac{4 \pi s}{5}) + 
\frac{2}{49}(\cos\frac{2 \pi s}{7}+\cos\frac{4 \pi s}{7}+\cos\frac{6 \pi s}{7}) 
,\nonumber \\
V(s,{\cal S}_8) & = & \frac{s^7}{203212800} -
\frac{17\cdot s^5}{9676800} + \frac{s^3}{8294400}\cdot (1343+225\cos \pi s)+
 \nonumber \\
&&s\cdot (-\frac{16133}{4976640}-\frac{1}{256} \cos \frac{\pi s}{2} + 
\frac{1}{243} \cos \frac{2\pi s}{3} -
\frac{31}{12288} \cos \pi s)+  \nonumber \\
&&\frac{1}{32}(\sin \frac{\pi s}{4} - \sin \frac{3\pi s}{4})-
\frac{1}{128} \sin \frac{\pi s}{2} 
+\frac{1}{162\sqrt{3}} \sin \frac{2 \pi s}{3} 
+\frac{1}{18\sqrt{3}} \sin \frac{\pi s}{3} + \nonumber \\
&&\frac{4}{125}(\sin \frac{2 \pi}{5}\sin \frac{4 \pi s}{5}-
\sin \frac{\pi}{5}\sin \frac{2 \pi s}{5}) - \nonumber \\
 &&
\frac{1}{49}(\sin \frac{2 \pi s}{7}\csc \frac{\pi}{7}-
\sin \frac{4 \pi s}{7}\csc \frac{2 \pi}{7}+
\sin \frac{6 \pi s}{7}\csc \frac{3 \pi}{7}) 
,\nonumber \\
V(s,{\cal S}_9) & = & \frac{s^8}{14631321600} -
\frac{19\cdot s^6}{418037760} + 
\frac{145597\cdot s^4}{16721510400}+
\frac{s^3}{73728}\cdot \sin \pi s -
 \nonumber \\
&&s^2\cdot (\frac{67293991}{140460687360} + 
\frac{1}{4374} \cos \frac{2\pi s}{3}) -\nonumber \\
&&s\cdot (\frac{1}{256\sqrt{2}} \sin \frac{\pi s}{2} + 
\frac{1}{1458\sqrt{3}} \sin \frac{2 \pi s}{3} +
\frac{205}{98304} \sin \pi s)+ \frac{199596951167}{56184274944000}+ \nonumber \\
&&
\frac{1}{64}(\cos \frac{\pi s}{4} \csc \frac{\pi}{8}
- \cos \frac{3\pi s}{4} \csc \frac{3 \pi}{8})+
\frac{2}{125} (\cos \frac{4 \pi s}{5}-
\cos \frac{2 \pi s}{5}) -\nonumber \\
&&
\frac{5}{512\sqrt{2}}\cos \frac{\pi s}{2} 
+\frac{257}{17496} \cos \frac{2 \pi s}{3} 
+\frac{1}{36\sqrt{3}} \cos \frac{\pi s}{3} + \nonumber \\
&&\frac{2}{81}
(-\cos \frac{2 \pi s}{9}+\cos \frac{4 \pi s}{9}+
 \cos \frac{8 \pi s}{9}) - \nonumber \\
 &&
\frac{1}{98}
(\cos \frac{2 \pi s}{7}
\csc \frac{\pi}{7}\csc \frac{2\pi}{7}+
\cos \frac{4 \pi s}{7}\csc \frac{2 \pi}{7}
\csc \frac{3 \pi}{7}+
\cos \frac{6 \pi s}{7}\csc \frac{3 \pi}{7}
\csc \frac{\pi}{7}) 
,\nonumber \\
V(s,{\cal S}_{10}) & = & \frac{s^9}{1316818944000} -\frac{11\cdot s^7}{12541132800}+
\frac{113113\cdot s^5}{358318080000}-\frac{\sin \pi s}{2949120}\cdot s^4  -\nonumber \\
&& \frac{18063859\cdot s^3}{468202291200}+s^2\cdot (\frac{1}{4374 \sqrt{3}}
\sin \frac{2 \pi s}{3} +\frac{143}{1179648} \sin \pi s)+\nonumber \\
&& s\cdot [\frac{273512277643}{240789749760000}+
\frac{1}{512 \sqrt{2}}\cos \frac{\pi s}{2}+\frac{7}{13122}\cos \frac{2 \pi s}{3}+ \nonumber \\
&& \frac{1}{625}(\cos \frac{4 \pi s}{5}-\cos \frac{2 \pi s}{5}) ]-
\frac{2877523}{707788800} \sin \pi s -\frac{1211}{52488 \sqrt{3}}\sin \frac{2 \pi s}{3}-\nonumber \\
&& \frac{5}{1024 \sqrt{2}} \sin \frac{\pi s}{2} -\frac{1}{108} \sin \frac{\pi s}{3}+
\frac{1}{64 \sqrt{2}}(\csc \frac{3 \pi}{8} \sin \frac{3 \pi s}{4}-
\csc \frac{\pi}{8} \sin \frac{\pi s}{4})+\nonumber \\
&& \frac{1}{50}(\sin \frac{3 \pi s}{5}-\sin \frac{\pi s}{5})-
\frac{2\sqrt{2}}{625} (\frac{\sqrt{5}+2}{\sqrt{5+\sqrt{5}}} \sin \frac{2 \pi s}{5}+
\frac{\sqrt{5}-2}{\sqrt{5-\sqrt{5}}} \sin \frac{4 \pi s}{5})-\nonumber \\
&& \frac{1}{196} \csc \frac{\pi}{7} \csc \frac{2 \pi}{7}\csc \frac{3 \pi}{7}
(\sin \frac{6 \pi s}{7}+\sin \frac{4 \pi s}{7}-\sin \frac{2 \pi s}{7})+\nonumber \\
&& \frac{1}{81} (\csc \frac{4 \pi}{9} \sin \frac{8 \pi s}{9}+
\csc \frac{2 \pi}{9} \sin \frac{4 \pi s}{9}+
\csc \frac{\pi}{9} \sin \frac{2 \pi s}{9}).
\nonumber
\end{eqnarray}

\subsection{Coxeter groups.}

Let us define two auxiliary functions 
\begin{eqnarray}
U_{+}(s,p,{\sl G})&=&V(s+p,{\sl G})+V(s-p,{\sl G})\;,\nonumber \\
U_{-}(s,p,{\sl G})&=&V(s+p,{\sl G})-V(s-p,{\sl G})\;
\label{new}
\end{eqnarray}
with obvious properties
\begin{eqnarray}
U_{+}(s,p,{\bf d}^m/d_r)&=&U_{-}(s,p+\frac{d_r}{2},{\bf d}^m)-U_{-}(s,p-\frac{d_r}{2},{\bf d}^m)\;,\;\;
U_{+}(s,0,{\sl G})=2 V(s,{\sl G}),\nonumber \\
U_{-}(s,p,{\bf d}^m/d_r)&=&U_{+}(s,p+\frac{d_r}{2},{\bf d}^m)-U_{+}(s,p-\frac{d_r}{2},{\bf d}^m)\;,\;\;
U_{-}(s,\frac{d_r}{2},{\bf d}^m)=V(s,{\bf d}^m/d_r),\nonumber
\label{new2}
\end{eqnarray}
where $(m-1)$-$tuple$ $\{{\bf d}^m/d_r\}=\{d_1,d_2,...,d_{r-1},d_{r+1},...,d_m\}$ doesn't contain $d_r$-exponent.

Sylvester waves for the Coxeter groups are given below expressed through the 
relations elaborated in the previous Sections.

$\underline {G={\sl A_m}}$ $\;,\;\;\;\;\;\;d_r=2,3,..., m+1\;;
\;\;\;\;\;\;\;\xi({\sl A_m})=\frac{1}{4}m(m+3)$
\begin{eqnarray}
V(s,{\sl A}_m)&=&U_{-}(s,\frac{1}{2},{\cal S}_m).
\label{coxaa}
\end{eqnarray}

$\underline {G={\sl B_m}}$ $\;,\;\;\;\;\;\;d_r=2,4,6,...,2 m\;;
\;\;\;\;\;\;\;\;\xi({\sl B_m})=\frac{1}{2}m(m+1)$
\begin{eqnarray}
V(s,{\sl B}_m)&=&\frac{1}{2}\Psi_2(s-\xi({\sl B_m}))\cdot 
U_{+}(\frac{s}{2},0,{\cal S}_m)\;.
\label{coxb1}
\end{eqnarray}

In the list for ${\sl D_m}$ groups the degree $m$ occurs twice when $m$ is even.
This is the only case involving such a repetition.

$\underline {G={\sl D_m}}$ $\;,\;\;\;\;\;\;d_r=2,4,6,...,2 (m-1),m\;,\;\;m \geq 3\;;
\;\;\;\;\;\;\xi({\sl D_m})=\frac{1}{2}m^2\;\;,\;\;$
\begin{eqnarray}
V(s,{\sl D}_{2m})&=&\Psi_2(s)\cdot U_{+}(\frac{s}{2},\frac{m}{2},
{\cal S}_{2m}),\\
V(s,{\sl D}_{2m+1})&=&\sum_{s_1=0}^{s-\xi({\sl D_{2m+1}})} V(s+\frac{2m+1}{2}-s_1,{\sl B}_{2m})
\cdot \Psi_{2m+1}(s_1),  \nonumber \\
V(s,{\sl D}_{3})&=&V(s,{\sl A}_{3}),  \nonumber \\
V(s,{\sl D}_{5})&=&U_{-}(s,\frac{11}{2},{\cal S}_8)-
U_{-}(s,\frac{9}{2},{\cal S}_8)-
U_{-}(s,\frac{5}{2},{\cal S}_8)+U_{-}(s,\frac{3}{2},{\cal S}_8).\nonumber
\label{coxd1}
\end{eqnarray}

$\underline {G={\sl G_2}}$ $\;,\;\;\;\;\;\;d_r=2,6\;;\;\;\;\;\;\;\;\;\;\;\;\;\;\;\;\;
\xi({\sl G_2})=4\;,$
\begin{eqnarray}
V(s,{\sl G}_2)&=&\Psi_2(s)\cdot U_{-}(\frac{s}{2},1,{\cal S}_3).
\label{g21}
\end{eqnarray}

$\underline {G={\sl F_4}}$ $\;,\;\;\;\;\;\;d_r=2,6,8,12\;;\;\;\;\;\;\;\;\;\;\xi({\sl F_4})=14\;,$
\begin{eqnarray}
V(s,{\sl F}_4)&=&\Psi_2(s) \cdot [\;U_{+}(\frac{s}{2},\frac{7}{2},{\cal S}_6)-
U_{+}(\frac{s}{2},\frac{3}{2},{\cal S}_6)\;].
\label{f41}
\end{eqnarray}

$\underline {G={\sl E_6}}$ $\;,\;\;\;\;\;\;d_r=2,5,6,8,9,12\;;\;\;\;\;\;\;\;\;\;\;\;\;\;
\xi({\sl E_6})=21\;,$
\begin{eqnarray}
V(s,{\sl E}_6)&=&U_{+}(s,18,{\cal S}_{12})-U_{+}(s,17,{\cal S}_{12})-U_{+}(s,15,{\cal S}_{12})+
\nonumber \\
&&U_{+}(s,13,{\cal S}_{12})+U_{+}(s,5,{\cal S}_{12})-U_{+}(s,2,{\cal S}_{12})\;.
\label{e6}
\end{eqnarray}

$\underline {G={\sl E_7}}$ $\;,\;\;\;\;\;\;d_r=2,6,8,10,12,14,18\;;\;\;\;\;\;\;
\;\;\;\xi({\sl E_7})=35\;,$
\begin{eqnarray}
V(s,{\sl E}_7)&=&\Psi_2(s-1)\cdot [\; U_{+}(\frac{s}{2},5,{\cal S}_9)-
U_{+}(\frac{s}{2},3,{\cal S}_9)\;]\;.
\label{e71}
\end{eqnarray}

$\underline {G={\sl E_8}}$ $\;,\;\;\;\;\;\;d_r=2,8,12,14,18,20,24,30\;;\;\;\;\;\;\;\xi({\sl E_8})=64\;,$

\begin{eqnarray}
V(s,{\sl E}_8)&=&\Psi_2(s)\cdot [\;U_{-}(\frac{s}{2},28,{\cal S}_{15})+
U_{-}(\frac{s}{2},21,{\cal S}_{15})+U_{-}(\frac{s}{2},12,{\cal S}_{15})+\nonumber \\
&&U_{-}(\frac{s}{2},11,{\cal S}_{15})-U_{-}(\frac{s}{2},8,{\cal S}_{15})-U_{-}(\frac{s}{2},7,
{\cal S}_{15})-\nonumber \\
&&U_{-}(\frac{s}{2},6,{\cal S}_{15})-U_{-}(\frac{s}{2},26,{\cal S}_{15})-U_{-}(\frac{s}{2},25,
{\cal S}_{15})\;].
\label{e81}
\end{eqnarray}

$\underline {G={\sl H_3}}$ $\;,\;\;\;\;\;\;d_r=2,6,10\;;\;\;\;\;\;\;\;\;\;
\;\;\;\xi({\sl H_3})=9\;,$
\begin{eqnarray}
V(s,{\sl H}_3)&=&\Psi_2(s-1)\cdot
[\;U_{+}(\frac{s}{2},3,{\cal S}_5)-U_{+}(\frac{s}{2},1,{\cal S}_5)\;].
\label{h31}
\end{eqnarray}

$\underline {G={\sl H_4}}$ $\;,\;\;\;\;\;\;d_r=2,12,20,30\;;\;\;\;
\;\;\;\xi({\sl H_3})=32\;,$
\begin{eqnarray}
V(s,{\sl H}_4)&=&U_{+}(s,32,{\sl E}_8)-U_{+}(s,24,{\sl E}_8)-U_{+}(s,18,{\sl E}_8)-
U_{+}(s,14,{\sl E}_8)+\nonumber \\
&&U_{+}(s,10,{\sl E}_8)-U_{+}(s,8,{\sl E}_8)+U_{+}(s,6,{\sl E}_8)+U_{+}(s,0,{\sl E}_8).
\label{h41}
\end{eqnarray}


$\underline {G={\sl I}_{m}}$ $\;,\;\;\;\;\;\;d_r=2,m\;;\;\;\;\;\;\;\;\;\;\xi({\sl I}_{m})=
1+\frac{1}{2}m$
\begin{eqnarray}
V(s,{\sl I_m})&=&\sum_{s_1=0}^{s-\xi({\sl I}_m)} 
\Psi_2(s-\xi({\sl I}_{m})-s_1)\cdot \Psi_m(s_1)\;,\\ \nonumber
V(s,{\sl I_2})&=&V(s,{\sl B_1})\;,\;\;
V(s,{\sl I_3})= V(s,{\sl A_2})\;,\;\;
V(s,{\sl I_4})= V(s,{\sl B_2})\;,\\ \nonumber
V(s,{\sl I_5})&=&U_{+}(s,\frac{7}{2},{\sl A_4})-U_{+}(s,\frac{1}{2},{\sl A_4})\;,\\ \nonumber
V(s,{\sl I_6})&=&V(s,{\sl G_2})\;,\;\;
V(s,{\sl I_8})= U_{+}(s,5,{\sl B_4})-U_{+}(s,1,{\sl B_4})\\ \nonumber
V(s,{\sl I_{10}})&=& U_{-}(s,3,{\sl H_3})\;,\;\;
V(s,{\sl I_{12}})= U_{+}(s,7,{\sl F_4})-U_{+}(s,1,{\sl  F_4})\;.
\label{im1}
\end{eqnarray}

\section{Acknowledgement}

We would like to thank Prof. Z.Rudnik for the communication of the proof related to the 
asymptotic behaviour of the {\it least common multiple} {\sf lcm}($1,2,...,m,N$).

This research was supported in part by grants from the Tel Aviv University Research
Authority and the Gileadi Fellowship program of the Ministry of Absorption
of the State of Israel (LGF).

\newpage



\newpage

\appendix
\renewcommand{\theequation}{\thesection\arabic{equation}}
\section{Asymptotic behaviour of ${\sf lcm}$(1,2,...,N).}
\label{appendix}
\setcounter{equation}{0}

The arithmetical function {\it least common multiple} {\sf lcm}($1,2,...,N$)
=${\cal L}(N)$ of the series of the natural numbers takes a specific place among the
other arithmetical
functions. It can neither be represented as the Cauchy integral of the generating function
with subsequent evaluation with Hardy-Ramanujan circle method like different partition
functions $p(N),\;q(N)$, nor has it its genesis in Riemann's Zeta-function like many
arithmetical functions $\mu(N),\;\nu(N),\;\phi(N),\;d(N)$.
${\cal L}(N)$ appears naturally in the theory of restricted partition numbers
as periods of Sylvester waves in symmetric groups ${\cal S}_N$ .

\begin{figure}[h]
\psfig{figure=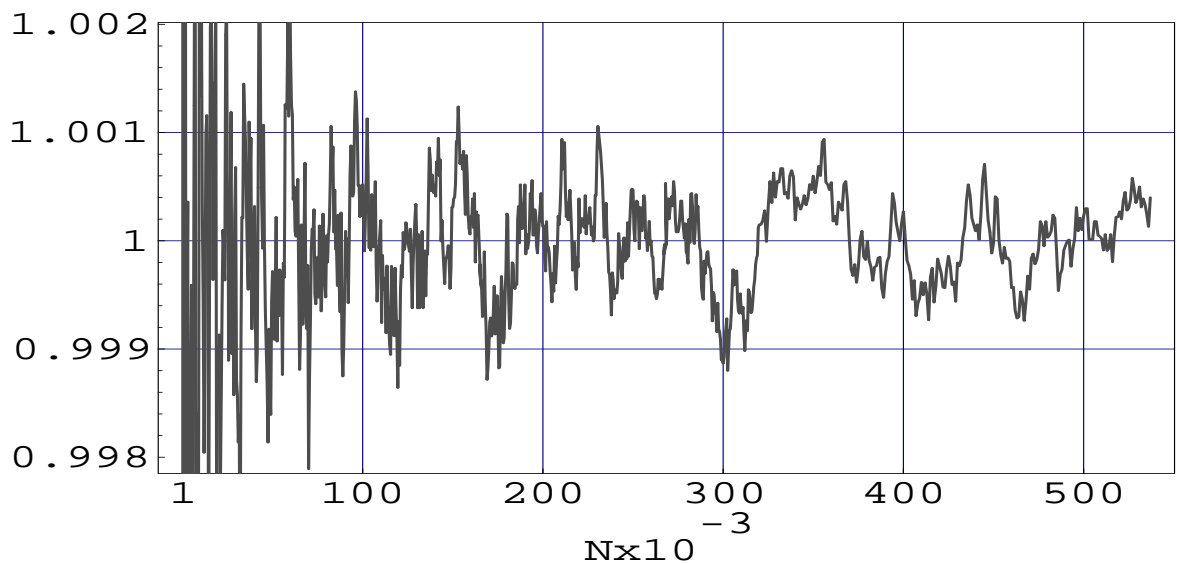,height=8cm,width=17cm}
\caption{Asymptotic behaviour of $\frac{1}{N}\ln[{\cal L}(N)]$.}
\label{lcm}
\end{figure}

Numerical calculations of $\frac{1}{N}\ln[{\cal L}(N)]$ in the range
$0<N<550\times 10^3$ give an oscillating behaviour around 1 with asymptotic approach to this value
(Fig. \ref{lcm}).
This enabled us to conjecture an asymptotic law
\begin{equation}
\lim_{N\to\infty} \frac{\ln {\cal L}(N)}{N}= 1.
\label{zeev1}
\end{equation}
In the rest of this Appendix we give a proof of this statement.
Before going to the proof we recall some facts of the prime number theory:

${\sf F1.}$ The Prime Number Theorem (PNT)
\begin{equation}
\mbox{if}\;\;\;\;
\pi(N) = \sum_{p_i \leq N} 1 \;\;\;\;,\;\;\;\;\;\;
\mbox{then}\;\;\;\pi(N)\stackrel{N\rightarrow \infty} \simeq \frac{N}{\ln N}\;.
\label{pnt1}
\end{equation}
where a sum is running over all primes $p_i$ up to $N$.

${\sf F2.}$ Let us set after Chebyshev
\begin{equation}
\theta (N) = \sum_{p_i \leq N} \ln p_i\;,
\label{pnt2}
\end{equation}
then PNT is equivalent to $\theta (N)\stackrel{N\rightarrow \infty} \simeq N$.

${\sf F3.}$ The Rieman hypothesis is equivalent to
\begin{equation}
\theta (N) = N+ O(\sqrt{N} \ln N)
\label{pnt3}
\end{equation}
Now it follows

${\sf Lemma \; A.}$
\[\lim_{N\to\infty} \frac{\ln {\cal L}(N)}{N}= 1.\]
$ and\;\; assuming\;\; the\;\; Rieman\;\; hypothesis$
\[\ln {\cal L}(N)= N+ O(\sqrt{N} \ln N).\]

${\sf Proof\;\;of\;\;Lemma\; A .}$

We write the prime decomposition of ${\cal L}(N)$ as ${\cal L}(N)=\prod p^{k_p}$.
Clearly , for a prime to divide ${\cal L}(N)$ it has to be at most less than $N$.
Moreover the highest $k$ power of $p$ dividing one of the integers 1, 2, ..., $N$ is
\[ k_p=[\frac{\ln N}{\ln p}] \]
Thus we find
\begin{equation}
\ln {\cal L}(N)=\sum_{p_i \leq N} [\frac{\ln N}{\ln p}]\cdot \ln p
\label{lem1}
\end{equation}
To estimate $\ln {\cal L}(N)$, break the previous sum into two parts, one $Q_1$
coming from primes $p \leq \sqrt{N}$ and the second $Q_2$ from primes $\sqrt{N}\leq p \leq N$ :
\[ Q_1=\sum_{p_i \leq \sqrt{N}} [\frac{\ln N}{\ln p}]\cdot \ln p\;,\;\;
Q_2=\sum_{\sqrt{N}\leq p \leq N} [\frac{\ln N}{\ln p}]\cdot \ln p\;\]
For estimating $Q_1$, use $[x]\leq x$ and so we find
\[ Q_1\leq \sum_{p_i \leq \sqrt{N}} \frac{\ln N}{\ln p}\cdot \ln p=
\ln N\cdot \pi(\sqrt{N})\simeq 2 \sqrt{N}\]
by the PNT. For the second sum $Q_2$ note that if $\sqrt{N}\leq p \leq N$ then
\[ 1\leq \frac{\ln N}{\ln p} < 2 \]
and hence its integer part is identically 1. Thus
\[ Q_2=\sum_{\sqrt{N}\leq p \leq N} 1\cdot \ln p=\theta (N)-\theta (\sqrt{N}) \]
Since $\theta (\sqrt{N}) \simeq \sqrt{N}$ , we obtain finally\[ \ln {\cal L}(N)= \theta (N)+\theta (\sqrt{N}) \]
Our ${\sf Lemma\; A}$ follows immediately from {\sf F2, F3}. $\;\;\;\;\;\;\blacksquare$

\newpage

\renewcommand{\theequation}{\thesection\arabic{equation}}
\section{Derivation of Sylvester waves $V(s,{\cal S}_4)$ and $V(s,{\cal S}_5)$.}
\label{appendix}
\setcounter{equation}{0}

We will illustrate how do the formulas (\ref{knf}-\ref{knf13}) work in the case of the 
symmetric groups ${\cal S}_4$ and ${\cal S}_5$. 

We start with Sylvester wave $V(s,{\cal S}_3)$ taken from
(\ref{sym4})
\begin{equation}
V(s,{\cal S}_3)=\frac{s^2}{12} - \frac{7}{72} -\frac{1}{8} \cos \pi s +
\frac{2}{9} \cos \frac{2 \pi s}{3}
\label{b1}
\end{equation}
and with successive usage of the formulas (\ref{knf}) and (\ref{knf11}) one can obtain
\begin{eqnarray}
R_1^4(s)=\frac{1}{144}\;,\;\;R_2^4(s)=0\;,\;\;R_3^4(s)=-\frac{1}{96}\cdot (5+3 \cos \pi s)\;,\;\;\;\;
{\cal R}_4^4(s)=-\frac{2}{9\sqrt{3}} \sin \frac{2 \pi s}{3}\;.
\label{b2}
\end{eqnarray}
Now we will use the representation (\ref{knf12})
\begin{equation}
V(s,{\cal S}_4)=\sum_{j=1}^{3}R^{4}_j(s)\cdot s^{4-j}+{\cal R}^{4}_{4}(s)+
\rho^4_1 \cdot \sin \frac{\pi}{2}s+\rho^4_2 \cdot \sin \pi s\;.
\label{b3}
\end{equation}
Since $V(s,{\cal S}_4)=W(s-5,{\cal S}_4)$ the variable $s$ takes only integer values 
what makes
the last contribution in (\ref{b3}) into the $V(s,{\cal S}_4)$ irrelevant. 
The unknown coefficient
$\rho^4_1$ is determined with help of {\it zeroes} (\ref{lem31}) of $W(s,{\cal S}_4)$
\begin{equation}
0=V(1,{\cal S}_4)=\sum_{j=1}^{3}R^{4}_j(1)+{\cal R}^{4}_{4}(1)+ \rho^4_1\;\;,\;\;\;\;\mbox{or}
\;\;\;\;\rho^4_1=\frac{1}{8}
\label{b4}
\end{equation}
Thus we arrive at the Sylvester wave $V(s,{\cal S}_4)$ presented in (\ref{sym4}).

Repeating the same procedure with symmetric group ${\cal S}_5$ we find
\begin{eqnarray}
R_1^5(s)&=&\frac{1}{2880}\;,\;\;R_2^5(s)=0\;,\;\;R_3^5(s)=-\frac{11}{1152}\;,\;\;R_4^5(s)=-\frac{1}{64}\sin \pi s\;,\\ 
{\cal R}_5^5(s)&=&\frac{475}{27648}-\frac{2}{27}\cos \frac{2 \pi s}{3}+\frac{1}{8\sqrt{2}}\cos \frac{\pi
s}{2}\;.\nonumber
\label{b6}
\end{eqnarray}
The representation (\ref{knf13}) produces
\begin{equation}
V(s,{\cal S}_5)=\sum_{j=1}^{4}R^{5}_j(s)\cdot s^{5-j}+{\cal R}^{5}_{5}(s)+
\rho^5_0 +\rho^5_1 \cdot \cos \frac{2 \pi s}{5}+\rho^5_2 \cdot \cos \frac{4 \pi s}{5}\;.
\label{b7}
\end{equation}
Since $V(s,{\cal S}_5)=W(s-\frac{15}{2},{\cal S}_5)$ the variable $s$ has only half-integer
values. By solving three linear equations $V(\frac{1}{2},{\cal S}_5)=V(\frac{3}{2},{\cal S}_5)=
V(\frac{5}{2},{\cal S}_5)=0$ we find
\begin{equation}
\rho^5_0=\frac{217}{28800}\;,\;\;\rho^5_1=-\frac{2}{25}\;,\;\;\rho^5_2=\frac{2}{25}\;,
\label{b8}
\end{equation}
which together with (\ref{b7}) produces the Sylvester wave 
$V(s,{\cal S}_5)$ from (\ref{sym4}).

\newpage

\section{Table of restricted partition numbers $W(s,{\cal S}_m)$.} 

In this Appendix we give the Table of the restricted partition numbers 
${\cal P}_m(s)=W(s,{\cal S}_m)\;m \leq 10$ for $s$ running in the different ranges. 
One can verify that the content of this Table can be obtained  with the help of the formulas
(\ref{sym4}).

\begin{center}
\begin{tabular}{|c|c|c|c|c|c|c|c|c|c|c|c|} \hline\hline
$s\;$ & & $S_1$ & $S_2$ & $S_3$ & $S_4$ & $S_5$ & $S_6$  & $S_7$  & $S_8$ & $S_9$ & $S_{10}$ \\ \hline\hline
 1   & &  1    &   1   &   1   &   1    &   1   &   1    &   1    &   1   &   1   &   1      \\ \hline
 2   & &  1    &   2   &   2   &   2    &   2   &   2    &   2    &   2   &   2   &   2      \\ \hline  
 3   & &  1    &   2   &   3   &   3    &   3   &   3    &   3    &   3   &   3   &   3      \\ \hline
 4   & &  1    &   3   &   4   &   5    &   5   &   5    &   5    &   5   &   5   &   5      \\ \hline
 5   & &  1    &   3   &   5   &   6    &   7   &   7    &   7    &   7   &   7   &   7      \\ \hline
 6   & &  1    &   4   &   7   &   9    &   10  &   11   &   11   &   11  &   11  &   11     \\ \hline
 7   & &  1    &   4   &   8   &  11    &   13  &   14   &   15   &   15  &   15  &   15     \\ \hline
 8   & &  1    &   5   &  10   &  15    &   18  &   20   &   21   &   22  &   22  &   22     \\ \hline
 9   & &  1    &   5   &  12   &  18    &   23  &   26   &   28   &   29  &   30  &   30     \\ \hline  
10   & &  1    &   6   &  14   &  23    &   30  &   35   &   38   &   40  &   41  &   42     \\ \hline\hline 
51   & &  1    &  26   &  243  &  1215  &  4033 &  9975  & 19928  & 33940 & 51294 & 70760    \\ \hline
52   & &  1    &  27   &  252  &  1285  &  4319 & 10829  & 21873  & 37638 & 57358 & 79725    \\ \hline
53   & &  1    &  27   &  261  &  1350  &  4616 & 11720  & 23961  & 41635 & 64015 & 89623    \\ \hline
54   & &  1    &  28   &  271  &  1425  &  4932 & 12692  & 26226  & 46031 & 71362 & 100654   \\ \hline
55   & &  1    &  28   &  280  &  1495  &  5260 & 13702  & 28652  & 50774 & 79403 & 112804   \\ \hline
56   & &  1    &  29   &  290  &  1575  &  5608 & 14800  & 31275  & 55974 & 88252 & 126299   \\ \hline
57   & &  1    &  29   &  300  &  1650  &  5969 & 15944  & 34082  & 61575 & 97922 & 141136   \\ \hline
58   & &  1    &  30   &  310  &  1735  &  6351 & 17180  & 37108  & 67696 & 108527 & 157564  \\ \hline
59   & &  1    &  30   &  320  &  1815  &  6747 & 18467  & 40340  & 74280 & 120092 & 175586  \\ \hline
60   & &  1    &  31   &  331  &  1906  &  7166 & 19858  & 43819  & 81457 & 132751 & 195491  \\ \hline\hline
101  & &  1    &   51  &  901  &  8262  & 48006 & 198230 & 628998 &1621248&3539452&6757864   \\ \hline
102  & &  1    &   52  &  919  &  8505  & 49806 & 207338 & 662708 &1719877&3778074&7254388   \\ \hline
103  & &  1    &   52  &  936  &  8739  & 51649 & 216705 & 697870 &1823402&4030512&7782608   \\ \hline
104  & &  1    &   53  &  954  &  8991  & 53550 & 226479 & 734609 &1932418&4297682&8345084   \\ \hline
105  & &  1    &   53  &  972  &  9234  & 55496 & 236534 & 772909 &2046761&4580087&8942920   \\ \hline
106  & &  1    &   54  &  990  &  9495  & 57501 & 247010 & 812893 &2167057&4878678&9578879   \\ \hline
107  & &  1    &   54  &  1008 &  9747  & 59553 & 257783 & 854546 &2293142&5194025&10254199  \\ \hline
108  & &  1    &   55  &  1027 &  10018 & 61667 & 269005 & 898003 &2425678&5527168&10971900  \\ \hline
109  & &  1    &   55  &  1045 &  10279 & 63829 & 280534 & 943242 &2564490&5878693&11733342  \\ \hline
110  & &  1    &   56  &  1064 &  10559 & 66055 & 292534 & 990404 &2710281&6249733&12541802  \\
\hline\hline
\end{tabular}
\end{center}

\end{document}